\documentclass{article}
\usepackage{amssymb,amsmath}
\usepackage{graphicx}

\def\qed{\hfill {\hbox{${\vcenter{\vbox{               %HOLLOW SQUARE
   \hrule height 0.4pt\hbox{\vrule width 0.4pt height 6pt
   \kern5pt\vrule width 0.4pt}\hrule height 0.4pt}}}$}}}

\newtheorem{theorem}{Theorem}
\newtheorem{definition}{Definition}

\newtheorem{proposition}[theorem]{Proposition}

\newtheorem{example}{Example}

\newtheorem{remark}{Remark}

{\qed\par\medskip}

\date{}

\title{\Large \textbf{Birack modules and their link invariants}}

\author{Regina Bauernschmidt\footnote{Email: rrb02007@pomona.edu} \and 
Sam Nelson\footnote{Email: knots@esotericka.org}}

\begin{document}
\maketitle

\begin{abstract}
We extend the rack algebra $\mathbb{Z}[X]$ defined by Andruskiewitsch and 
Gra\~na to the case of biracks, enabling
a notion of birack modules. We use these birack modules to define an
enhancement of the birack counting invariant generalizing the birack module
counting invariant in \cite{HHNYZ}. We provide examples demonstrating
that the enhanced invariant is not determined by the Jones or Alexander 
polynomials and is strictly stronger than the unenhanced birack counting
invariant.
\end{abstract}

\medskip

\quad
\parbox{5in}{
\smallskip

}

\section{\large\textbf{Introduction}}

\textit{Biracks} were first introduced in \cite{FRS} as a generalization of
\textit{racks}, algebraic structures whose axioms encode blackboard-framed
isotopy of oriented knots and links in $\mathbb{R}^3$. Initially called 
``wracks'' alluding to the ``wrack and ruin'' resulting from keeping only
conjugation in a group, racks were first considered by Conway and Wraith
in the 1950s \cite{FR}. The current term ``rack'' (without the ``w'') is due to
Fenn and Rourke in \cite{FR}, who dropped the ``w'' to denote dropping 
invariance under writhe-changing type I Reidemeister moves.

In \cite{N2}, 
an integer-valued invariant of knots and links associated to a finite birack 
$X$ known as the \textit{integral birack counting invariant}, 
$\Phi_X^{\mathbb{Z}}$ was introduced by the second author, generalizing the 
rack counting invariant from \cite{N}. In \cite{AG} an associative algebra 
known as the
\textit{rack algebra} $\mathbb{Z}[X]$ was defined for every finite rack, and
in \cite{HHNYZ} the rack counting invariant from \cite{N} was enhanced
with representations (known as \textit{rack modules}) of a modified form 
of $\mathbb{Z}[X]$.

In this paper, for every finite birack $X$ we define an associative algebra 
$\mathbb{Z}_B[X]$ we call the \textit{birack algebra}. We use representations 
of $\mathbb{Z}_B[X]$, known as \text{birack modules} or $X$-modules, to 
enhance the birack counting invariant from \cite{N2}.
The new invariant is defined for classical and virtual knots and links.

The paper is organized as follows: in section \ref{B} we review the basics 
of biracks and the birack counting invariant. In section \ref{BA} we 
define the birack algebra and give examples of birack modules. In section
\ref{I} we use birack modules to enhance the birack counting invariant
and give examples to demonstrate that the birack module enhanced invariant
is not determined by the Jones or Alexander polynomials and is strictly 
stronger than the unenhanced birack counting invariant. We conclude with
a few questions for future research in section \ref{Q}.

\section{\large\textbf{Biracks}}\label{B}

We begin with a definition (see \cite{FJK, N2}).

\begin{definition}\textup{
Let $X$ be a set. A \textit{birack structure} on $X$ is an invertible
map $B:X\times X\to X\times X$ satisfying
\begin{list}{}{}
\item[(i)] $B$ is \textit{sideways invertible}, that is, there exists a 
unique invertible map $S:X\times X\to X\times X$ satisfying for all $x,y\in X$
\[S(B_1(x,y),x)=(B_2(x,y),y),\]
\item[(ii)] $B$ is \textit{diagonally invertible}, that is, the compositions 
$S^{\pm 1}_{1}\circ \Delta$ and $S^{\pm 1}_{2}\circ \Delta$
of  the diagonal map $\Delta:X\to X\times X$ given by $\Delta(x)=(x,x)$ 
with the components of the sideways map and its inverse are bijections, and
\item[(iii)] $B$ is a solution to the \textit{set-theoretic Yang-Baxter 
equation:}
\[(B\times I)(I\times B)(B\times I)=
(I\times B)(B\times I)(I\times B).\]
\end{list}
We will often abbreviate $B_1(x,y)=y^x$ and $B_2(x,y)=x_y$.
}\end{definition}

\begin{example}
\textup{Perhaps the simplest example of a birack structure is the
class of \textit{constant action biracks}: let $X$ be a set and let 
$\sigma,\tau:X\to X$ be bijections. Then the map $B(x,y)=(\sigma(y),\tau(x))$
defines a birack structure on $X$ iff $\sigma$ and $\tau$ commute; see 
\cite{N2}.}
\end{example}

\begin{example}
\textup{Let $\tilde\Lambda=\mathbb{Z}[t^{\pm 1},s,r^{\pm 1}]/(s^2-s(1-tr))$. Then
any $\tilde\Lambda$-module $M$ is a birack, known as $(t,s,r)$\textit{-birack},
under the map $B(x,y)=(ty+sx,rx)$. This structure generalizes the \textit{$(t,s)$-rack} structure from \cite{FR}. The idea is to consider the case when the
output components $B$ are linear combinations of the inputs, i.e., set 
$B(x,y)=(sx+ty,rx+py)$; then the invertibility conditions (i) and (ii) are 
satisfied if one of the four coefficients must be zero. Selecting $p=0$, 
the Yang-Baxter equation is then satisfied provided we have $s^2=(1-tr)s$.}
\end{example}

It is frequently useful to specify a birack structure on a finite set
$X=\{x_1,x_2,\dots, x_n\}$ by listing the operation tables for the components
of $B$ viewed as binary operations $(x,y)\mapsto B_1(y,x)=x^y$ and
$(x,y)\mapsto B_2(x,y)=x_y$. Note the reversed order of the operands in 
$B_1$; this is for compatibility with previous work representing these 
operations as right actions. Then a birack operation $B$ on $X$ is specified
by the $n\times 2n$ matrix $M_B=[U|L]$ with $U(i,j)=k$ and $L(i,j)=l$
where $x_k=B_1(x_j,x_i)$ and $x_l=B_2(x_i,x_j)$. Conversely, such a matrix
specifies a birack structure iff the operation it defines via
\[B(x_i,x_j)=(x_{U[j,i]},x_{L[i,j]})\] 
satisfies the birack axioms.

\begin{example}\textup{
The constant action birack on $X=\{x_1,x_2,x_3\}$ with bijections 
$\sigma=(123)$ and $\tau=(132)$ has birack matrix
\[M_B=\left[\begin{array}{ccc|ccc}
2 & 2 & 2 & 3 & 3 & 3 \\
3 & 3 & 3 & 1 & 1 & 1 \\
1 & 1 & 1 & 2 & 2 & 2 \\
\end{array}\right]\]
and the $(t,s,r)$-birack on $X=\mathbb{Z}_4$ with $t=r=3$ and $s=2$ has 
birack matrix
\[M_B=\left[\begin{array}{cccc|cccc}
1 & 3 & 1 & 3 & 1 & 1 & 1 & 1 \\
4 & 2 & 4 & 2 & 4 & 4 & 4 & 4 \\
3 & 1 & 3 & 1 & 3 & 3 & 3 & 3 \\
2 & 4 & 2 & 4 & 2 & 2 & 2 & 2 \\
\end{array}\right]\]
where $0=x_1, 1=x_2, 2=x_3,$ and $3=x_4$, using 4 as the representative of 
the class of 0 in $\mathbb{Z}_4$ since our row and column numbering starts 
with 1. Then for instance the entry in row 1 column 3 of the left matrix 
should be $B_1(x_3,x_1)=(tx_1+sx_3)=3(1)+2(3)=3+6=9=1$, while the entry in
row 3 column 2 of the right matrix is $B_2(x_3,x_2)=r(x_3)=2(2)=4$.
}\end{example}

As with other algebraic structures, we have the following standard notions:
\begin{definition}\textup{
Let $X$ and $X'$ be sets with birack operations $B:X\times X\to X\times X$
and $B':X'\times X'\to X'\times X'$. Then:
\begin{itemize}
\item a map $f:X\to X'$ is a \textit{birack homomorphism} if for all
$x,y\in X$ we have
\[B'(f(x),f(y))=(f(B_1(x,y)),f(B_2(x,y)));\]
\item a subset $Y\subset X$ is a \textit{subbirack} of $X$ if the restriction
of $B$ to $Y\times Y$ defines a birack structure on $Y$; equivalently,
$Y\subset X$ is a subbirack if $Y\times Y$ is closed under $B$.
\end{itemize}
}\end{definition}

The birack axioms encode the blackboard-framed Reidemeister moves
with semiarcs labeled according to the rules below.
\[\includegraphics{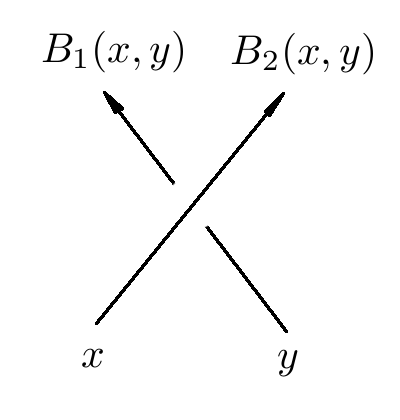} \quad
\includegraphics{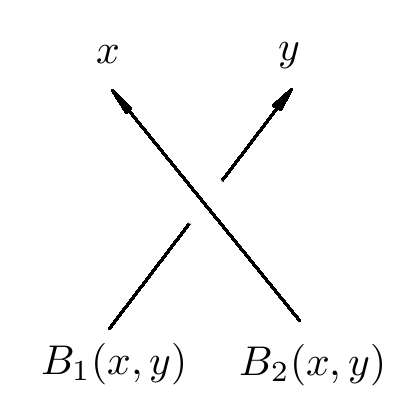} \quad
\includegraphics{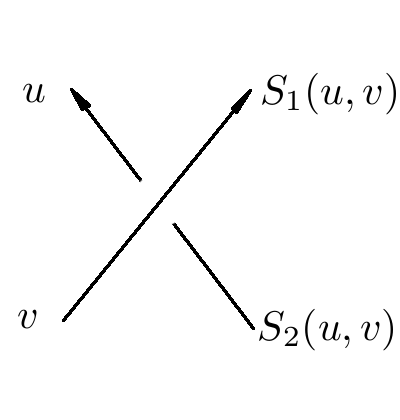} \quad
\includegraphics{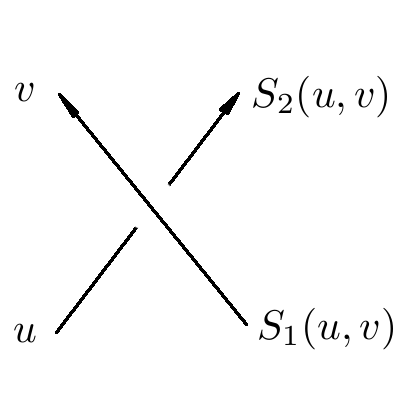}
\]
%See \cite{N} for more.

At a positive kink, the invertibility properties of $B$ define bijections
$\alpha:X\to X$ and $\pi: X\to X$ defined by 
$\alpha=(S_2^{-1}\circ \Delta)^{-1}$ and $\pi=S_1^{-1}\circ\Delta\circ \alpha$
which give us the labels on the other semiarcs in a blackboard framed
type I move.
\[\includegraphics{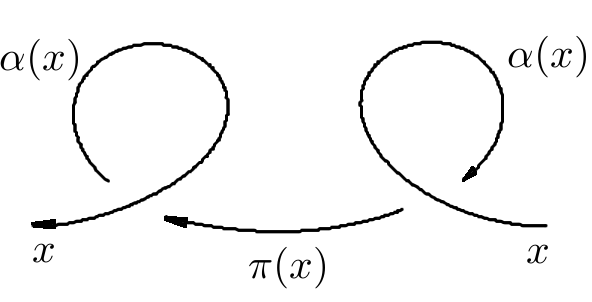}\]
The exponent $N$ of the bijection $\pi$, that is, the smallest positive
integer $N$ satisfying $\pi^N(x)=x$ for all $x\in X$, is known as the
\textit{birack rank} or \textit{birack characteristic} of $X$. If $X$ is
a finite set, then $N$ is guaranteed to exist, but infinite biracks may
have infinite rank.

\begin{example}
\textup{A birack with rank $N=1$ is called a \textit{strong biquandle} (see 
\cite{FJK, KR, NV}); a birack satisfying $B_2(x,y)=x$ is a \textit{rack}
(see \cite{FR}). A strong biquandle which is also a rack is a \textit{quandle};
see \cite{J,M}.}
\end{example}

\begin{example}
\textup{Every oriented blackboard framed link diagram $L$ has a 
\textit{fundamental birack} $FB(L)$, defined as the set of equivalence 
classes of birack words in the semiarcs of $L$ modulo the equivalence relation
determined by the birack axioms and the crossing relations in $L$. For 
example, the knot $5_1$ below has the fundamental birack presentation listed
below:}
\[\includegraphics{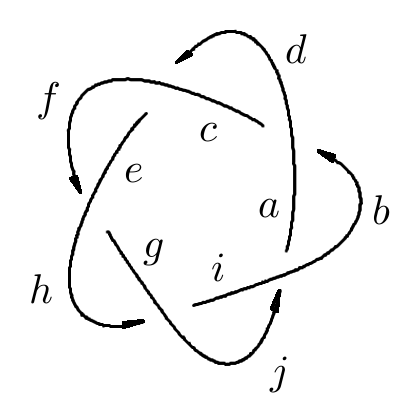}\raisebox{0.75in}{
$\begin{array}{rl} 
FB(5_1)=\langle a,\ b,\ c,\ d,\ e,\ f,\ g,\ h,\ i,\ j\ | & 
B(a,b)=(c,d),\ B(c,d)=(e,f), \\ 
& B(e,f)=(g,h), \ B(g,h)=(i,h),\\ & B(i,j)=(a,b)\rangle.\end{array}$}\]
\end{example}

A labeling of semiarcs in a framed oriented link diagram with elements
of a birack $X$ determines a homomorphism $f:FB(L)\to X$ if and only
if the birack labeling condition is satisfied by the labels in $X$ at every
crossing. Conversely, every homomorphism $f:FB(L)\to X$ determines a 
unique labeling of the semiarcs of $L$. By construction, we have the 
following standard result:

\begin{theorem}
If $L$ and $L'$ are oriented link diagrams which are related by blackboard
framed Reidemeister moves and $X$ is a finite birack, then the sets of
labelings of the semiarcs of $L$ and $L'$ by elements of $X$ satisfying
the birack labeling condition, $\mathrm{Hom}(FB(L),X)$ and 
$\mathrm{Hom}(FB(L'),X)$, are in bijective correspondence.
\end{theorem}

Thus, the cardinality of the set of birack labelings of an oriented link 
diagram by a finite birack $X$ is a positive integer valued invariant of
framed isotopy. To obtain an invariant of unframed ambient isotopy, we
observe as in \cite{N} that the cardinalities of the sets of birack 
labelings are periodic in the birack rank $N(X)$ since birack labelings
are preserved by the \textit{$N$-phone cord move}:
\[\includegraphics{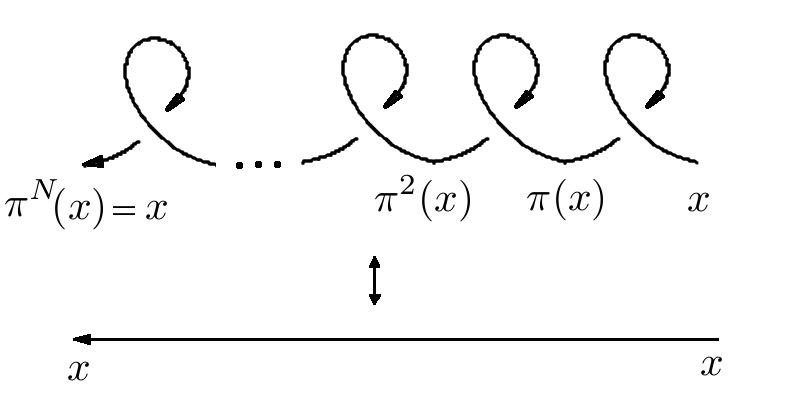}\]

Thus, summing these
cardinalities over a complete period of framings mod $N$ yields an 
invariant of ambient isotopy known as the \textit{integral birack counting
invariant} $\Phi_X^{\mathbb{Z}}.$ More formally, we have

\begin{definition}
\textup{Let $X$ be a finite birack with birack rank $N$ and $L$ an oriented
link with $c$ components. Let $W=(\mathbb{Z}_N)^c$ and let $(L,\mathbf{w})$
be a blackboard framed diagram of $L$ with writhe $w_k$ on component $k$ for
each $k=1,\dots,c$. The \textit{integral birack counting invariant} is the 
integer}
\[\Phi_X^{\mathbb{Z}}(L)
=\sum_{\mathbf{w}\in W} |\mathrm{Hom}(FB((L,\mathbf{w})),X)|.\]
\end{definition}

\begin{example}\label{ex1}\textup{
Let $X$ be the birack structure on $X=\{1,2\}$ with birack matrix
\[M_B=\left[\begin{array}{cc|cc}
2 & 2 & 2 & 2 \\
1 & 1 & 1 & 1 \\
\end{array}\right].\]
As a labeling rule, this birack requires switching labels from $1$ to $2$
and $2$ to $1$ on both semiarcs at every crossing point. Thus, a component 
in a diagram has 
two valid labelings iff it has an even number of semiarcs, but since a 
classical diagram with $n$ crossings has $2n$ semiarcs, every diagram has 
2 valid labelings. Since $N=1$, we have $\Phi_X^{\mathbb{Z}}(L)=2^c$ for
every classical link with $c$ components.}
\[\includegraphics{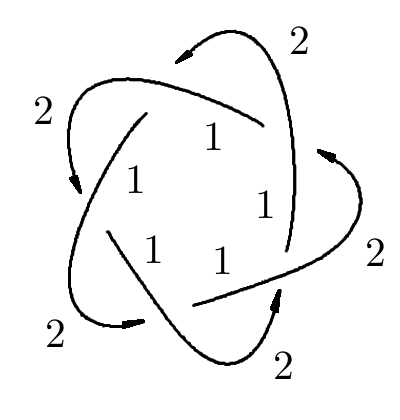} \quad \includegraphics{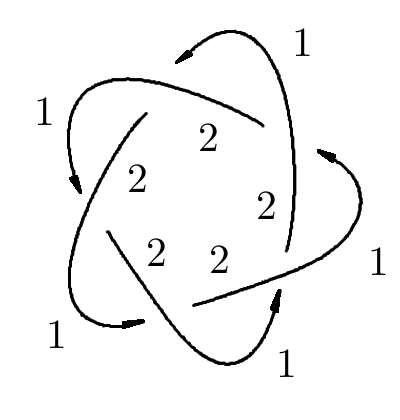} 
\quad \raisebox{0.7in}{$\Phi_{X}^{\mathbb{Z}}(5_1)=2$}\]
\[\scalebox{0.95}{
\includegraphics{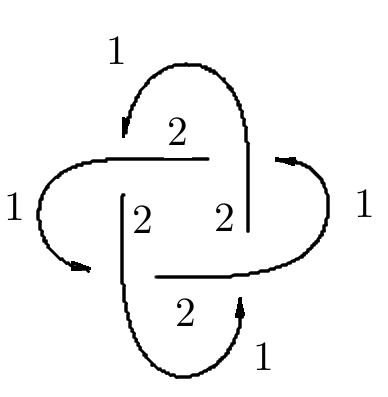} \quad \includegraphics{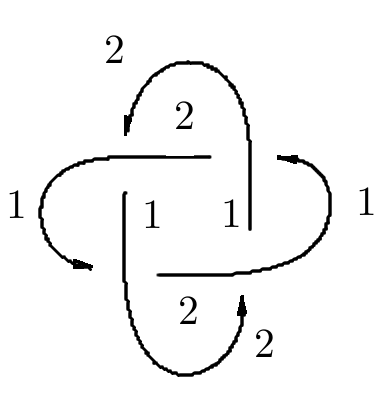} 
\includegraphics{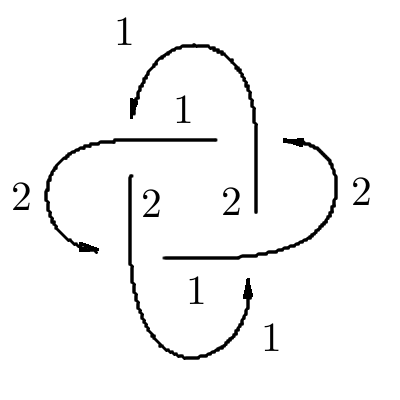} \quad \includegraphics{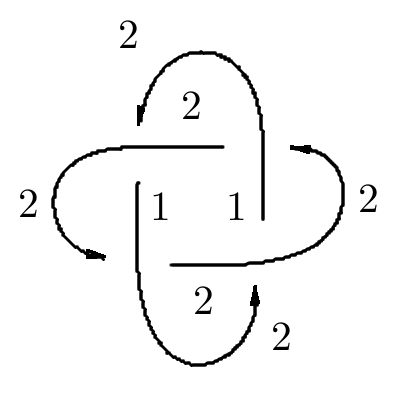} }
\quad \raisebox{0.7in}{$\Phi_{X}^{\mathbb{Z}}(L4a1)=4$}\]
\end{example}

The integral birack counting invariant is also defined for virtual links
via the usual technique of ignoring virtual crossings when dividing the link
into semiarcs or, equivalently, regarding our virtual link diagrams as
being drawn on a surface with sufficient genus to avoid virtual crossings.

\begin{example}
\textup{While the birack in example \ref{ex1} does not distinguish classical
links with the same number of components via $\Phi_X^{\mathbb{Z}}$, it 
\textit{can} distinguish some virtual links from other virtual links with 
the same number of components. Recall that a \textit{virtual link diagram}
can include virtual crossings represented by circled self-intersections which
we may regard as indicating genus in the surface on which the link diagram is
drawn.
For instance, the \textit{virtual Hopf link}
$\mathrm{vHopf}$ below has no valid $X$-labelings, so 
$\Phi_X^{\mathbb{Z}}(\mathrm{vHopf})=0\ne 2=
\Phi_X^{\mathbb{Z}}(\mathrm{Unlink})$.}
\[\includegraphics{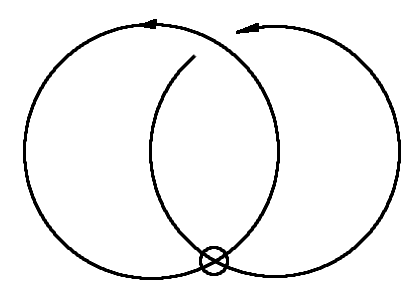}\]
\textup{See \cite{K} for more about virtual knots and links.}
\end{example}

\section{\large\textbf{Birack algebras and modules}}\label{BA}

We would like to enhance the integral birack counting invariant by finding
a way to distinguish between labelings rather than just counting how many
total labeling we have. To this end, we will use a scheme analogous to the
$(t,s,r)$-birack structure to define an associative algebra we call the
\textit{birack algebra} associated to a finite birack $X$. We will discuss 
the motivation for this definition later in this section.

\begin{definition}\textup{Let $X$ be a finite birack of birack rank $N$. 
The \textit{birack
algebra} associated to $X$ is the quotient of the polynomial algebra
$\mathbb{Z}[t^{\pm 1}_{x,y}, s_{x,y}, r^{\pm 1}_{x,y}]$ generated by
invertible elements $t_{x,y},r_{x.y}$ and general elements $s_{x,y}$
for all $x,y\in X$, modulo the ideal $I$ defined by the relations}
\begin{itemize}
\item{$r_{x_y,z}r_{x,y}-r_{x_{z^y},y_z}r_{x,z^y}$}
\item{$t_{x_{z^y},y_z}r_{y,z}-r_{y^x,z^{x_y}}t_{x,y}$}
\item{$s_{x_{z^y},y_z}r_{x,z^y}-r_{y^x,z^{x_y}}s_{x,y}$}
\item{$t_{x,z^y}t_{y,z}-t_{y^x,z^{x_y}}t_{x_y,z}$}
\item{$t_{x,z^y}s_{y,z}-s_{y^x,z^{x_y}}t_{x,y}$}
\item{$s_{x,z^y}-t_{y^x,z^{x_y}}s_{x_y,z}r_{x,y}-s_{y^x,z^{x_y}}s_{x,y}$}
\item{$1-\displaystyle{\prod_{k=0}^{N-1} 
(t_{\pi^k(x),\alpha(\pi^k(x))}r_{\pi^k(x),\alpha(\pi^k(x))}+s_{\pi^k(x),\alpha(\pi^k(x))})}$}
\end{itemize}
\textup{We will denote the birack algebra of $X$ by $\mathbb{Z}_B[X]$.
A representation of $\mathbb{Z}_B[X]$, i.e. an abelian group $G$ with
automorphisms $t_{x,y},r_{x,y}:G\to G$ and endomorphisms $s_{x,y}:G\to G$
such that the above listed maps are zero, will be called a 
\textit{birack module} or simply an $X$-module. }
\end{definition}

\begin{remark}\textup{
This definition generalizes the quandle and rack module definitions from 
\cite{AG}, which were used in \cite{CEGS,HHNYZ} to enhance the quandle 
and rack counting invariants respectively. In particular, if $X$ is a rack
then a rack module in the sense of \cite{HHNYZ} is a birack module over
$X$ in which every $r_{x,y}=\mathrm{Id}_G$. Note that in general, the set 
of rack modules
over a rack $X$ considered as a birack is a proper subset of the set of 
birack modules over $X$.}\end{remark}

The motivation behind the birack algebra definition is to define secondary
labelings of $X$-labeled oriented blackboard framed link diagrams with
a ``bead'' at every semiarc and to use the $(t,s,r)$-birack operations on the
beads at a crossing with each $t,s$ and $r$ coefficient indexed by the
$X$-labels on the input strands at a positive crossing and the output
strands at a negative crossing.
\[\includegraphics{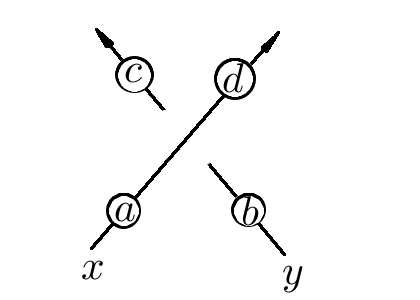} \quad \includegraphics{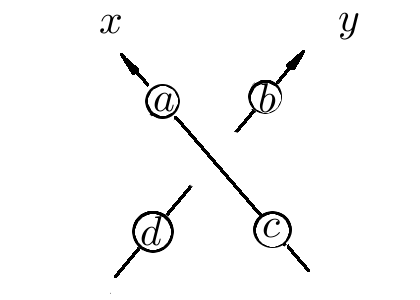} \quad 
\raisebox{0.5in}{$\begin{array}{rcl}
c & = & s_{x,y} a+t_{x,y} b \\
d & = & r_{x,y} a \\
\end{array}$}
\]

The birack algebra relations are chosen to preserve bead labelings under
the blackboard-framed Reidemeister moves and the $N$-phone cord move. The 
choice of bead labeling rules guarantees that for every blackboard framed 
$X$-labeled diagram, the number of bead labelings is the same before and 
after type II moves and framed type I moves:
\[\includegraphics{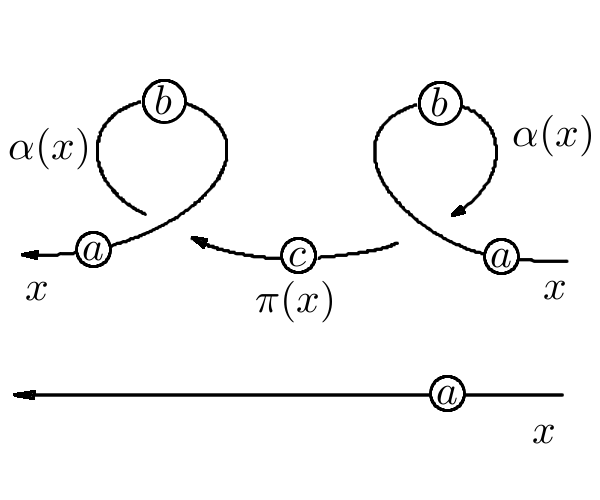} \quad \raisebox{0.9in}{
$\begin{array}{rcl}
b & = & r_{x,\alpha(x)} a \\
c & = & s_{x,\alpha(x)} a+t_{x,\alpha(x)} b \\
\Rightarrow c & = & (s_{x,\alpha(x)}+t_{x,\alpha(x)}r_{x,\alpha(x)})a
\end{array}$}
\]
\[\includegraphics{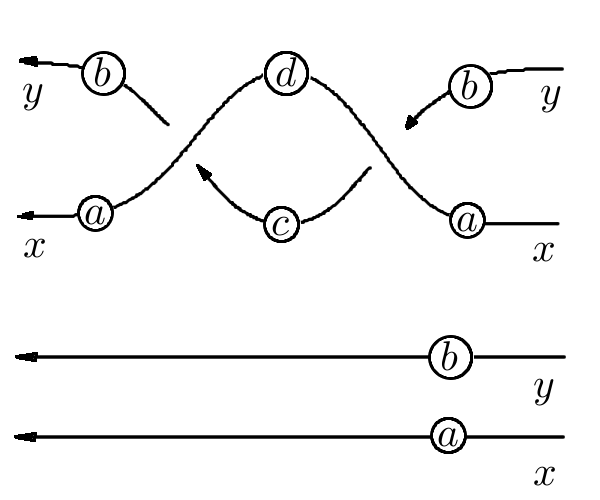} \quad \raisebox{0.7in}{
$\begin{array}{rcl}
c & = & s_{x,y} a+t_{x,y} b \\
d & = & r_{x,y} a \\
\end{array}$}\]
The other type II and framed type I cases are similar.

Six of the seven birack algebra relations come from the type III move:
\[\includegraphics{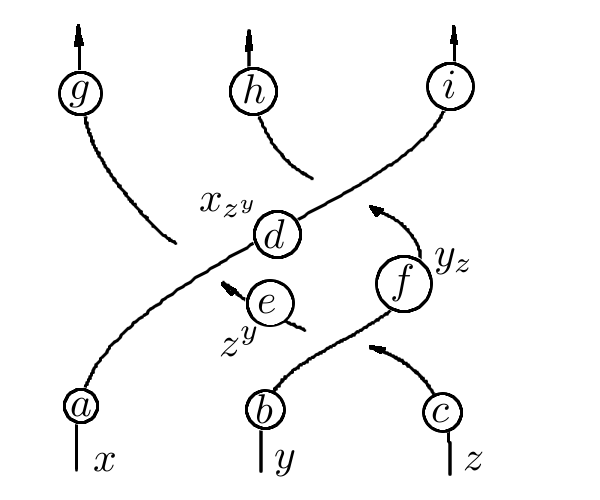}\quad \includegraphics{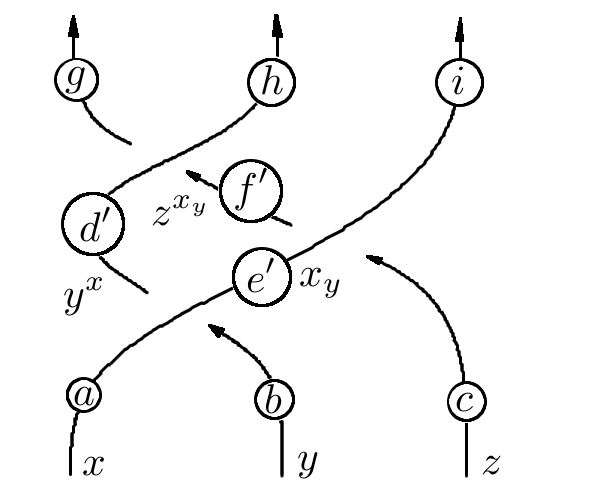}\]
Comparing coefficients on the output beads $g,h$ and $i$ yields the 
first six relations.

The final birack algebra relation comes from the $N$-phone cord move.
Pushing a bead on a strand labeled with birack element $x$ though a 
positive kink multiplies the bead by 
$(t_{x,\alpha(x)}r_{x,\alpha(x)}+s_{x,\alpha(x)})$; thus, we need the 
product of these over a complete period of framings mod $N$ to be $1$.

\begin{example}\textup{
Let $X$ be a finite birack and $R$ a ring. We can give $R$ the structure 
of a $\mathbb{Z}_B[X]$-module by choosing invertible elements 
$t_{x,y},r_{x,y}\in R$ and elements $s_{x,y}$ for all $x,y\in X$ such that
the birack algebra relations are satisfied. We can represent such a 
birack module structure with a \textit{birack module matrix} $M_R=[T|S|R]$
where $T(i,j)=t_{x_i,x_j}$, etc. For example, the birack from example 
\ref{ex1} has birack modules on $\mathbb{Z}_3$ including}
\[M_R=\left[\begin{array}{rr|rr|rr}
2 & 2 & 0 & 2 & 2 & 1 \\
2 & 2 & 2 & 0 & 1 & 2 \\
\end{array}\right].\]
\end{example}

\begin{example}\textup{
Another important example of a $\mathbb{Z}_B[X]$-module is the 
\textit{fundamental $\mathbb{Z}_B[X]$-module} of an $X$-labeled oriented
blackboard framed link diagram $L$, denoted $\mathbb{Z}_f[X]$, where
$f:FB(L)\to X$ is the $X$-labeling of $L$. Starting with our $X$-labeled
diagram, we put a bead on every semiarc and obtain a system of linear
equations with coefficients in $\mathbb{Z}_B[X]$ determined at the crossings.
For convenience, we will represent such a module with the coefficient matrix 
of the homogeneous system. For example, the link $L4a1$ from example
\ref{ex1} with the labeling by the birack from the same example below
has the listed fundamental $\mathbb{Z}_B[X]$-module.}
\[\includegraphics{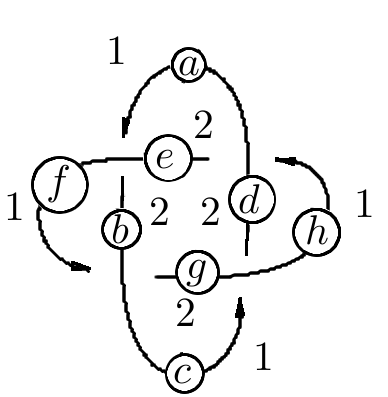}
\quad\raisebox{0.5in}{$\mathbb{Z}_f[X]=
\left[\begin{array}{cccccccc}
t_{2,1} & -1 & 0 & 0 & s_{2,1} & 0 & 0 & 0 \\
0 & 0 & 0 & 0 & r_{2,1} & -1 & 0 & 0 \\
0 & s_{2,1} & 0 & 0 & 0 & t_{2,1} & 0 & -1 \\
0 & r_{2,1} & -1 & 0 & 0 & 0 & 0 & 0 \\
0 & 0 & t_{2,1} & -1 & 0 & 0 & s_{2,1} & 0 \\
0 & 0 & 0 & 0 & 0 & 0 & r_{2,1} & -1 \\
0 & 0 & 0 & s_{2,1} & -1 & 0 & 0 & t_{2,1} \\
-1 & 0 & 0 & r_{2,1} & 0 & 0 & 0 & 0 \\
\end{array}\right]$}\]
\textup{The columns in this matrix represent the beads on the semiarcs, 
while each row represents an equation in the beads determined at a crossing. 
For example, the upper left crossing determines the equations 
$t_{2,1}a+s_{2,1}e=b$ and $r_{2,1}e=f$, which give us the first two rows of the 
matrix.} 
\end{example}

\section{\large\textbf{Birack module enhancements of the counting invariant}}\label{I}

We will now use an $X$-module $R$ to enhance the integral birack
counting invariant by taking the set  of $\mathbb{Z}_B[X]$-module 
homomorphisms $\mathrm{Hom}_{\mathbb{Z}[X]}(\mathbb{Z}_f[X],R)$
as a signature for each birack labeling $f\in\mathrm{Hom}(FB(L),X)$. 
Since each homomorphism in $\mathrm{Hom}_{\mathbb{Z}[X]}(\mathbb{Z}_f[X],R)$
can be understood as a labeling of the beads by elements of $R$, we 
are effectively counting bead labelings of $X$-labelings of $L$, enhancing 
the birack counting invariant with a bead counting invariant of $X$-labeling
of $L$.

\begin{definition}\textup{
Let $X$ be a finite birack of rank $N$, $R$ a $\mathbb{Z}[X]$-module, and
$L$ an oriented (classical or virtual) link of $c$ components. Let 
$W=(\mathbb{Z}_N)^c$ and for each $\mathbf{w}\in W$, let $(L,\mathbf{w})$ 
be a diagram of $L$ with writhe $w_k$ on component $k$. Then the 
\textit{birack module enhanced multiset} is the multiset of $R$-modules}
\[\Phi_{X,R}^{B,M}(L)=\{\mathrm{Hom}_{\mathbb{Z}_B[X]}(\mathbb{Z}_f[X],R)\ |\ 
f\in \mathrm{Hom}(FB((L,\mathbf{w})),X), \mathbf{w}\in W\}\]
\textup{and the \textit{birack module enhanced polynomial} is}
\[\Phi_{X,R}^B(L)
=\sum_{\mathbf{w}\in W}\left(
\sum_{f\in\mathrm{Hom}(FB(L,\mathbf{w}),X)}u^{|\mathrm{Hom}_{\mathbb{Z}_B[X]}(\mathbb{Z}_f[X],R)|}\right).\]
\textup{If $R$ is not a finite set, we replace the cardinality 
$|\mathrm{Hom}_{\mathbb{Z}_B[X]}(\mathbb{Z}_f[X],R)|$ of the set of bead 
labelings with the rank of the 
$R$-module $\mathrm{Hom}_{\mathbb{Z}_B[X]}(\mathbb{Z}_f[X],R)$.}
\end{definition}

By construction, we have
\begin{proposition}
If $L$ and $L'$ are ambient isotopic oriented classical or virtual links,
then $\Phi_{X,R}^B(L)=\Phi_{X,R}^B(L')$ and 
$\Phi_{X,R}^{B,M}(L)=\Phi_{X,R}^{B,M}(L')$.
\end{proposition}

Note that the integral birack counting invariant $\Phi_{x}^{\mathbb{Z}}(L)$
is recovered from the birack module enhanced invariant 
$\Phi_{X,R}^B(L)$ by evaluating at $u=1$ and from $\Phi_{X,R}^{B,M}(L)$ by 
taking cardinality. In our next example we show that $\Phi_{X,R}^{B}(L)$
is stronger in general than $\Phi_{X}^{\mathbb{Z}}(L)$.

\begin{example}
\textup{Let $X=\mathbb{Z}_3$. We can give $X$ the structure of a 
$(t,s,r)$-birack by setting $t=1$, $s=2$ and $r=2$ so we have
$s^2=2^2=1=(1-2(1))2=(1-tr)s$. Then $tr+s=2+2=1$, so $X$ is a biquandle, and
we have biquandle matrix}
\[M_X=\left[\begin{array}{ccc|ccc}
1 & 3 & 2 & 1 & 1 & 1 \\
2 & 1 & 3 & 3 & 3 & 3 \\
3 & 2 & 1 & 2 & 2 & 2 
\end{array}\right]\]
\textup{where $x_1=0, x_2=1$ and $x_3=2$.
Our \texttt{python} computations reveal 320 $\mathbb{Z}_B[X]$-module 
structures on $R=\mathbb{Z}_5$, including for instance}
\[M_R=\left[\begin{array}{ccc|ccc|ccc}
4 & 4 & 4 & 3 & 1 & 2 & 2 & 2 & 2 \\
1 & 2 & 2 & 1 & 4 & 3 & 4 & 4 & 4 \\
3 & 3 & 1 & 4 & 3 & 2 & 1 & 1 & 1 \\
\end{array}\right]\]
\textup{The unknot has three $X$-labelings, one for each element of $X$,
as is easily seen by considering the zero-crossing unknot diagram. It is 
also easy to see that for each $X$-labeling, the space of bead-labelings
by $R$ is one-dimensional, with a total of 5 bead labelings for each
$X$-labeling. Thus, we have $\Phi_{X,R}^{B}(\mathrm{Unknot})=3u^5$.}

\textup{Now consider the virtual knot numbered $3.1$ in the \textit{Knot Atlas} 
\cite{KA}; it is the closure of the virtual braid diagram listed below. There 
are three classical crossings, all positive, and the bead and birack labels
do not change at the virtual crossings. Thus, we have six semiarcs and six 
equations for each $X$-labeling. There are three $X$-labelings of $3.1$, as
listed in the table. Recall that in our matrix notation from Example 3
we have $x_1=0, x_2=1$ and $x_3=2$, so the first row gives the monochomatic 
$X$-labeling with all zeroes, while the two other $X$-labelings are nontrivial.
For each $X$-labeling, we
replace the coefficients in the matrix of $\mathbb{Z}_f[X]$ with their values
from $M_R$ and row-reduce to find the contribution to $\Phi_{X,R}^{B}(3.1)$.}
\[\includegraphics{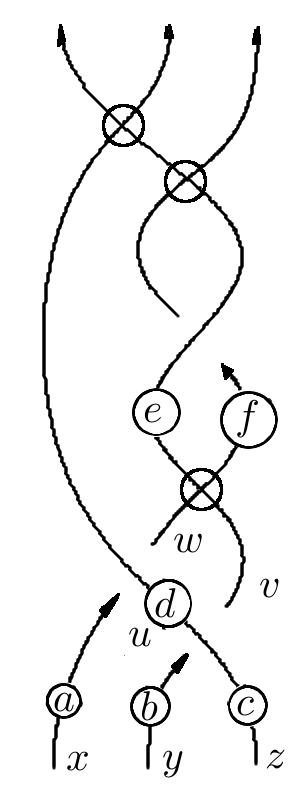} \quad\raisebox{1.5in}{$\begin{array}{c}
X-\mathrm{labelings:}\quad \begin{array}{|cccccc|}
\hline
x & y & z & u & v & w \\ \hline
1 & 1 & 1 & 1 & 1 & 1 \\
1 & 3 & 3 & 2 & 1 & 3 \\
1 & 2 & 2 & 3 & 1 & 3 \\ \hline
\end{array} \\
\ \\
M_{\mathbb{Z}_f[X]}=\left[\begin{array}{cccccc}
0 & -1 & 0 & s_{u,v} & t_{u,v} & 0 \\
0 & 0 & -1 & r_{u,v} & 0 & 0 \\
-1 & 0 & s_{z,w} & 0 & 0 & t_{z,w} \\
0 & 0 & r_{z,w} & -1 & 0 & 0 \\
0 & 0 & -1 & 0 & s_{v,w} & t_{v,w} \\
-1 & 0 & 0 & 0 & r_{v,w} & 0 \\
\end{array}\right] \\
\end{array}$}\]
\textup{The labeling of all semiarcs by $x_1\in X$ has matrix}
\[\left[\begin{array}{cccccc}
0 & -1 & 0 & s_{1,1} & t_{1,1} & 0 \\
0 & 0 & -1 & r_{1,1} & 0 & 0 \\
-1 & 0 & s_{1,1} & 0 & 0 & t_{1,1} \\
0 & 0 & r_{1,1} & -1 & 0 & 0 \\
0 & 0 & -1 & 0 & s_{1,1} & t_{1,1} \\
-1 & 0 & 0 & 0 & r_{1,1} & 0 \\
\end{array}\right]
\rightarrow \left[\begin{array}{cccccc}
0 & 4 & 0 & 3 & 4 & 0 \\
0 & 0 & 4 & 2 & 0 & 0 \\
4 & 0 & 3 & 0 & 0 & 4 \\
0 & 0 & 2 & 4 & 0 & 0 \\
0 & 0 & 4 & 0 & 3 & 4 \\
4 & 0 & 0 & 0 & 2 & 0 \\
\end{array}\right]
\rightarrow \left[\begin{array}{cccccc}
1 & 0 & 0 & 0 & 0 & 0 \\
0 & 1 & 0 & 0 & 0 & 0 \\
0 & 0 & 1 & 0 & 0 & 0 \\
0 & 0 & 0 & 1 & 0 & 0 \\
0 & 0 & 0 & 0 & 1 & 0 \\
0 & 0 & 0 & 0 & 0 & 1 \\
\end{array}\right]
\]
\textup{so there is only the zero-bead labeling and
this $X$-labeling contributes $u^1$ to $\Phi_{X,R}^{B}(3.1)$. Similar 
computations reveal that the other $X$-labelings also contribute $u^1$, so 
we have $\Phi_{X,R}^{B}(3.1)=3u$. Comparing this to the unknot, we have
$\Phi_{X,R}^{B}(\mathrm{Unknot})=3u^5\ne 3u=\Phi_{X,R}^{B}(3.1)$, so this
 example shows that $\Phi_{X,R}^{B}$ detects the knottedness of $3.1$. In 
fact, this example shows much more: (1) since $\Phi_X^{\mathbb{Z}}(3.1)=3=
\Phi_X^{\mathbb{Z}}(\mathrm{Unknot})$, this example shows that the birack 
module-enhanced invariant $\Phi_{X,R}^{B}$ is stronger than unenhanced
counting invariant $\Phi_X^{\mathbb{Z}}$, and (2) since $3.1$ and the 
unknot both have Jones polynomial equal to 1, this example shows that 
$\Phi_{X,R}^B$ is not determined by the Jones polynomial.}
\end{example}

\begin{example}
\textup{As we saw in example \ref{ex1}, the integral birack counting 
invariant with respect to the birack $X$ with matrix}
\[M_X=\left[\begin{array}{cc|cc} 
2 & 2 & 2 & 2 \\ 
1 & 1 & 1 & 1
\end{array}\right]\]
\textup{has value $\Phi_{X}^{\mathbb{Z}}(K)=2$ for all classical knots $k$.
However, the enhanced invariant $\Phi_{X,R}^{B}$ is effective at distinguishing
knots and links. We randomly selected a $\mathbb{Z}[X]$-module structure
on $\mathbb{Z}_3$ given by matrix}
\[M_R=\left[\begin{array}{cc|cc|cc}
2 & 2 & 0 & 2 & 1 & 1 \\
1 & 1 & 2 & 0 & 2 & 2 \\
\end{array}\right]\]
\textup{and used our \texttt{python} code to compute $\Phi_{X,R}^B$
for all prime knots with up to 8 crossings and all prime links with up to
7 crossings; the results are collected in the table below.}
\[\begin{array}{r|l}
\Phi_{X,R}^{B} & L \\ \hline
2u^3 & \mathrm{Unknot}, 4_1, 5_1, 5_2, 6_2, 6_3, 7_1, 7_2, 7_3, 7_5, 7_6, 8_1, 8_2, 8_3, 8_4, 8_6, 8_7, 8_8, 8_9, 8_{12}, 8_{13}, 8_{14}, 8_{16}, 8_{17} \\
2u^9 & 3_1, 6_1, 7_4, 7_7, 8_5, 8_{10}, 8_{11}, 8_{15}, 8_{19}, 8_{20}, 8_{21} \\
2u^{27} & 8_{18} \\
2u^3+2u^9 & L2a1, L4a1, L5a1, L6a2, L7a4, L7a6 \\
2u^3+2u^{27} & L7a2, L7a3, L7n1, L7n2 \\
4u^9 & L6a1, L6a3, L7a1, L7a5 \\
2u^3+6u^{27} & L6a4 \\ 
2u^3+6u^9 & L6n1, L7a7 \\
8u^9 & L6a5 \\
\end{array}\] 
\end{example}

\begin{example}\textup{
For our final example, we note that our \texttt{python} computations
also show that for the same birack $X$ from example \ref{ex1} with 
the $\mathbb{Z}[X]$-module structure on $\mathbb{Z}_3$ given by the
birack module matrix}
\[M_R=\left[\begin{array}{cc|cc|cc}
1 & 2 & 0 & 2 & 1 & 1 \\
2 & 1 & 2 & 0 & 1 & 1 \\
\end{array}\right]\]
\textup{$\Phi_{X,R}^{B}$ distinguishes the knots $8_{18}$ and $9_{24}$ which
both have Alexander polynomial 
\[\Delta(8_{18})=-t^3+5t^2-10t+13-10t^{-1}+5t^{-2}-t^{-3}=\Delta(9_{24})\]
while $\Phi_{X,R}^{B}(8_{18})=2u^{27}\ne 2u^9=\Phi_{X,R}^{B}(9_{24})$. In 
particular, $\Phi_{X,R}^B$ is not determined by the Alexander polynomial.
We note that this example also implies that $\Phi_{X,R}^B$ is not determined
by the generalized Alexander polynomial for virtual knots, since both 
$8_{18}$ and $9_{24}$ have generalized Alexander polynomial 0.}
\end{example}

\section{\large\textbf{Questions}}\label{Q}

We end with a few open questions for future research.

Are there examples of biracks $X$ and birack modules $R$ whose 
$\Phi_{x,R}^{B}$ invariants detect mutation or Homflypt-equivalent 
knots? Can $\Phi_{x,R}^{B}$ distinguish knots and links with
the same Khovanov or Knot Floer homology?
What, if any, is the relationship between birack module invariants
and Vassiliev invariants?

In \cite{TC}, biquandle labelings are extended to define counting invariants
of knotted surfaces in $\mathbb{R}^4$, and in particular it is shown that
surface biquandles are just biquandles. Does the same hold for the 
\textit{surface biquandle algebra} obtained by requiring bead-labeling
invariance under Roseman moves?

\bigskip

\noindent
\textsc{Department of Mathematical Sciences \\
Claremont McKenna College \\
850 Columbia Ave. \\
Claremont, CA 91711}

\end{document}